\documentclass[proceedings,submission%
]{dmtcs}

\usepackage[latin1]{inputenc}
\usepackage{amsmath,amssymb,latexsym,graphics,color,graphicx}
\usepackage{tabmac}
\usepackage[all]{xy}

\newlength{\cellsize}
\newtheorem{theorem}{Theorem}[section]

\newtheorem{definition}[theorem]{Definition}
\newtheorem{example}[theorem]{Example}

\newcommand{\cyc}{\mathrm{cyc}}
\newcommand{\height}{\mathrm{height}}
\newcommand{\id}{\mathrm{id}}
\newcommand{\Partition}{\mathcal{P}}
\newcommand{\sign}{\mathrm{sign}}
\newcommand{\Rec}{\mathcal{R}}
\newcommand{\shift}{\mathrm{shift}}
\newcommand{\Tab}{\mathcal{T}}
\newcommand{\weight}{\mathrm{weight}}
\newcommand{\word}{\mathrm{word}}
\newcommand{\et}{\tilde{e}}
\newcommand{\ft}{\tilde{f}}
\newcommand{\st}{\tilde{s}}

\begin{document}

\title[Fusion coefficients]{A combinatorial formula for fusion coefficients}

\author{Jennifer Morse\addressmark{1}
\thanks{Partially supported by the NSF grants DMS--0652641, DMS--0638625, DMS--1001898.}
\and 
Anne Schilling\addressmark{2}
\thanks{Partially supported by the NSF grants DMS--0652641, DMS--0652652, DMS--1001256.}}

\address{\addressmark{1}Department of Mathematics, Drexel University, Philadelphia, PA 19104, U.S.A\\
\addressmark{2}Department of Mathematics, University of California, One Shields
Avenue, Davis, CA 95616-8633, U.S.A.}


\maketitle

\begin{abstract}
\paragraph{Abstract.}
Using the expansion of the inverse of the Kostka matrix in terms of tabloids 
as presented by E{\u{g}}ecio{\u{g}}lu and Remmel, we show that the fusion 
coefficients can be expressed as an alternating sum over cylindric tableaux. 
Cylindric tableaux are skew tableaux with a certain cyclic symmetry.
When the skew shape of the tableau has a cutting point, meaning that the 
cylindric skew shape is not connected, or if its weight has at most two parts,
we give a positive combinatorial formula 
for the fusion coefficients. The proof uses a slight modification of a 
sign-reversing involution introduced by Remmel and Shimozono.  
We discuss how this approach may work in general.

\paragraph{R\'esum\'e.}
En utilisant l'expansion de l'inverse de la matrice Kostka en termes de tablo\"ides 
introduite par E{\u{g}}ecio{\u{g}}lu et Remmel, nous montrons que les coefficients de 
fusion peuvent \^etre exprim\'es comme une somme altern\'ee sur les tableaux cylindriques.
Les tableaux cylindriques sont des tableaux qui pr\'esentent une certaine sym\'etrie cyclique.
Lorsque la forme du tableau a un point de coupure, ce qui signifie que
la forme cylindrique n'est pas connect\'e, ou lorsque son poids a au
plus deux parts, nous donnons une formule combinatoire positive des
coefficients de fusion. La d\'emonstration utilise une l\'eg\`ere modification 
de l'involution qui change le signe introduite par Remmel et Shimozono.
Nous discutons comment cette approche pourrait fonctionner en g\'en\'eral.
\end{abstract}

\section{Introduction}

The famous \textit{Littlewood--Richardson} rule~\cite{LR:1934} provides a combinatorial expression for the coefficients
$c_{\lambda \mu}^\nu$ in the expansion of a product of Schur functions
\begin{equation} \label{equation.LR}
	s_\lambda s_\mu = \sum_\nu c_{\lambda \mu}^\nu s_\nu.
\end{equation}
It states that $c_{\lambda\mu}^\nu$ is equal to the number of column-strict tableaux of skew shape $\nu/\lambda$
and content $\mu$ whose column reading word is lattice. 
Here $\lambda,\mu,\nu$ are partitions and a column-strict tableau
of shape $\nu/\lambda$ is a filling of the skew shape which is weakly increasing across rows and strictly increasing
across columns. The content of a tableau or word is $\mu=(\mu_1,\mu_2,\ldots)$, where $\mu_i$ counts the number of $i$
in the tableau or word. Furthermore, a word is lattice if all right subwords have partition content.

In this paper, we consider the analogous problem for \textit{fusion coefficients}, which first appeared in the literature as the
structure constants of the Verlinde fusion algebra for the $\widehat{\mathfrak{sl}}_n$ Wess--Zumino--Witten models of 
level $\ell$~\cite{TUY:1989,Verlinde:1988}. Kac and Walton~\cite{Kac:1990,Walton:1990,Walton:1990a} provided an 
efficient algorithm for computing fusion coefficients for any type. In type 
$A_{n-1}$ of level $\ell$, their formula is expressed 
as an alternating sum of the Littlewood--Richardson coefficients
\begin{equation} \label{equation.Kac Walton}
	c_{\lambda \mu}^{\nu,\ell,n} = \sum_{\sigma \in \widehat{S}_n} \sign(\sigma) \; c_{\lambda \mu}^{\sigma(\nu+\rho)-\rho}\;,
\end{equation}
where the sum is over the affine symmetric group $\widehat{S}_n$ generated by $\langle \sigma_0,\sigma_1,\ldots,
\sigma_{n-1} \rangle$, $\rho=(n-1,n-2,\ldots,1,0)$, the symmetric group acts 
on compositions by permuting their entries, and 
$\sigma_0(\lambda)=(\lambda_n+\ell+n,\lambda_2,\ldots,\lambda_{n-1},\lambda_1-\ell-n)$.

A notorious problem has been to find a direct positive or combinatorial formula 
for the fusion coefficients as opposed to an alternating expression 
as in~\eqref{equation.Kac Walton}.  Many attempts have been made along these
lines.  Tudose~\cite{Tudose:2002} in her thesis gave a 
combinatorial interpretation when $\lambda$ or $\mu$ has at most two columns. 
For $n=2,3,$ positive formulas are known~\cite{BMW:1992} as well as when 
$\lambda$ and $\mu$ are rectangles~\cite{SS:2001}.
Knutson formulated a conjecture for the quantum Littlewood--Richardson coefficients as presented in~\cite{BKT:2003}
in terms of puzzles~\cite{KnutsonTao:2003}. It is known that the quantum cohomology structure coefficients
are related to fusion coefficients~\cite{A:1995,BCF:1999}. Coskun~\cite{Coskun:2009} gave a positive geometric rule to 
compute the structure constants of the cohomology ring of two-step flag varieties in terms of Mondrian tableaux. 

There are many other interpretations and appearances of fusion coefficients. 
For example, Goodman and Wenzl~\cite{GW:1990} showed that the fusion coefficients 
are related to the structure coefficients 
of Hecke algebras at roots of unity. This was used in~\cite{LM:2008}
to show that they are special cases of the structure coefficients of the $k$-Schur functions.
As mentioned earlier, they are also related to the quantum cohomology structure coefficients~\cite{A:1995,BCF:1999} and
intertwiners in vertex operator algebras~\cite{Wassermann:1998}.
Postnikov~\cite{Postnikov:2005} formulated the quantum cohomology ring in terms of the affine nilTemperley--Lieb
algebra, and Korff and Stroppel~\cite{Korff:2011a,KorffStroppel:2010} provided an analogous construction of the fusion
ring in terms of the affine local plactic algebra. 

The main result of this paper is a simple proof of a combinatorial formula for 
the fusion coefficients for a general class of partitions, which includes all previously known cases
(that is, the two-column case of Tudose~\cite{Tudose:2002}, $n=2,3$ of~\cite{BMW:1992}, and cases considered
by Postnikov~\cite{Postnikov:2005}). The proof uses 
the fusion Pieri rule to obtain an expression of the fusion coefficients in 
terms of cylindric tableaux (see Sections~\ref{section.tabloids} 
and~\ref{section.cyclic}).  We then amend a sign-reversing involution of Remmel 
and Shimozono~\cite{RS:1998} to cancel all negative terms (see 
Section~\ref{section.remmel}).  We finish with a discussion of how
this method could lead to a formula for fusion coefficients in
general.

To state the main result, several definitions are needed.
A partition $\lambda$ is of rank $n$ if it has at most $n$ parts. It 
is of level $\ell$ if $\lambda_1-\lambda_n \le \ell$.  Let 
$\Partition^{\ell,n}$ be the set of level $\ell$ and rank $n$ partitions. 
In addition, let $\Rec^{\ell,n}$ be the set of partitions of rank $n$ with
first part not exceeding $\ell$, that is, the set of partitions 
contained in a rectangle of width $\ell$ and height $n$.

\begin{theorem} \label{theorem.fusion}
Let $\ell,n$ be two positive integers, $\lambda,\mu \in \Rec^{\ell,n}$ and $\nu \in \Partition^{\ell,n}$ such that
$|\lambda|+|\mu| = |\nu|$.
\begin{enumerate}
\item \label{thm.cutting}
If $\nu/\lambda+\shift_{-\ell,n}(\nu/\lambda)$ is not connected 
(see Section~\ref{section.cyclic} for the definition of $\shift$),
then the level $\ell$ fusion coefficients $c_{\lambda \mu}^{\nu,\ell,n}$ for type $A_{n-1}$
are given in terms of usual Littlewood--Richardson coefficients for transformed partitions
\begin{equation} \label{equation.fusion LR}
	c_{\lambda \mu}^{\nu,\ell,n} = c_{\tilde{\lambda} \mu}^{\tilde{\nu}},
\end{equation}
with $\tilde{\lambda}$ and $\tilde{\nu}$ as in Definition~\ref{definition.cut}.
\item \label{thm.sl2}
If $\mu$ has at most two parts, then~\eqref{equation.fusion LR} holds with $\tilde{\lambda}$ and
$\tilde{\nu}$ as in Definition~\ref{definition.sl2}.
\end{enumerate}
\end{theorem}

The fusion coefficients enjoy many symmetries. For example:
\begin{itemize}
\item {\bf Columns of height $n$:} If $\lambda$ contains a column of height $n$, then
$c_{\lambda \mu}^{\nu,\ell,n} = c_{\tilde{\lambda} \mu}^{\tilde{\nu},\ell,n}$ where $\tilde{\lambda}$
(resp. $\tilde{\nu}$) is obtained from $\lambda$ (resp. $\nu$) by removing or adding a column of height $n$.
\item {\bf Level-rank duality:} Denoting by $\lambda^t$ the transpose partition of $\lambda$, we have
$c_{\lambda \mu}^{\nu,\ell,n} = c_{\lambda^t \mu^t}^{\nu^t,n,\ell}$. Here $\nu^t$ should be identified
with its cyclic analogue of attaching all parts $\nu^t_i$ for $i>\ell$ to $\nu^t_{i-\ell}$.
\item {\bf Strange duality:} For $\lambda \in \Rec^{\ell,n}$ denote by $\lambda^\vee$ the 
complement of $\lambda$ in the rectangle of size $\ell \times n$. Then the fusion coefficient labeled by the
complement partitions is related to $c_{\lambda \mu}^{\nu,\ell,n}$. This is best described using toric shapes,
see~\cite{Postnikov:2005}.
\item {\bf $S_3$ symmetry:}
The fusion coefficients $c_{\lambda \mu}^{\nu,\ell,n}$ are symmetric with respect to any permutation
of the partitions $\lambda,\mu,\nu$ (up to certain transformations on $\nu$ to put all partitions on the same footing). 
For more details see~\cite{Postnikov:2005}.
\end{itemize}
If under any of the above symmetries, one of the cases of Theorem~\ref{theorem.fusion} holds, a combinatorial
formula for the corresponding fusion coefficient follows. In particular, Theorem~\ref{theorem.fusion}~\eqref{thm.sl2}
under the level-rank duality is equivalent to the case when $\mu$ has at most two columns. This corresponds
to the case studied by Tudose~\cite{Tudose:2002}. We would like to point out that the proof given here (see 
Section~\ref{subsection.proof}) is much
simpler than the proof in~\cite{Tudose:2002} which involves many case checks.

\subsection*{Acknowledgments}
We thank Jason Bandlow for collaboration at the beginning of this project.
AS would also like to thank Catharina Stroppel for discussions during the 
Special Trimester ``On the Interaction of Representation Theory with Geometry and Combinatorics" at the Hausdorff Institut 
in Bonn, January - April 2011. Thanks also to Christian Korff, Chris Manon, Alexander Postnikov, and Peter Tingley for
helpful discussions, and Anders Buch and Nicolas Thi\'ery for their help in getting {\tt lrcalc}
integrated into {\tt Sage}~\cite{sage,sagecombinat}, which helped with computer experiments.

\section{Tabloids}
\label{section.tabloids}

To derive a formula for the fusion coefficients, we will use the well-known relation between Schur functions
$s_\mu$ and the homogeneous symmetric functions $h_\alpha$ (see for example~\cite[\S I.6, Table 1]{Macdonald:1995})
\begin{equation}
\label{equation.sh}
	s_\mu = \sum_\alpha K_{\alpha \mu}^{-1} \;h_\alpha,
\end{equation}
where $K_{\alpha \mu}^{-1}$ is the inverse of the Kostka matrix.
E{\u{g}}ecio{\u{g}}lu and Remmel~\cite{ERemmel:1990} gave an
interpretation for the entries in the inverse Kostka 
matrix using a combinatorial structure called tabloids. 
Note that these tabloids are different from the ones used in the representation
theory of the symmetric group.

The definition of a tabloid is a filling of a partition $\mu$
with certain shapes called ribbons.  A \textit{ribbon} is a connected 
skew shape which does not contain any $2\times 2$ squares.
The \textit{height} of a ribbon is one less than the number of occupied rows.
A \textit{tabloid} of shape $\mu$ is then
a tiling of $\mu$ by ribbons such that each ribbon 
contains a cell in the first column.
The \textit{weight} of a tabloid is $\beta=(\beta_1,\beta_2,\ldots)$, where $\beta_i$ is the length of the ribbon
starting in the $i$-th cell from the bottom in the first column of $\mu$. Here we use French notation for the shape $\mu$ placing 
the longest part of $\mu$ at the bottom. The \textit{sign} of a tabloid $T$ is $(-1)^{\height(T)}$, where the height of $T$ is the sum 
of the heights of all ribbons it contains.
The \textit{type} of a tabloid $T$ with weight $\beta$
is the partition $\alpha$ obtained by 
rearranging $\beta$ into non-increasing order.

\begin{example}{\rm
The four tabloids $T$ of shape $\mu = (3,2,1)$ are
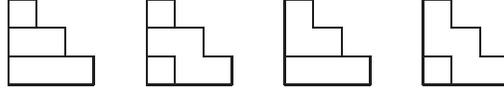
\begin{equation*}
	\unitlength=\cellsize
	\begin{picture}(3,3)
	\put(0,0){\line(1,0){3}}
	\put(0,1){\line(1,0){3}}
	\put(3,0){\line(0,1){1}}
	\put(0,2){\line(1,0){2}}
	\put(2,1){\line(0,1){1}}
	\put(0,3){\line(1,0){1}}
	\put(1,2){\line(0,1){1}}
	\put(0,0){\line(0,1){3}}
	\end{picture}%
	\qquad 
	\begin{picture}(3,3)
	\put(0,0){\line(1,0){3}}
	\put(0,1){\line(1,0){1}}
	\put(1,0){\line(0,1){1}}
	\put(3,0){\line(0,1){1}}
	\put(0,2){\line(1,0){2}}
	\put(2,1){\line(0,1){1}}
	\put(2,1){\line(1,0){1}}
	\put(0,3){\line(1,0){1}}
	\put(1,2){\line(0,1){1}}
	\put(0,0){\line(0,1){3}}
	\end{picture}%
	\qquad
	\begin{picture}(3,3)
	\put(0,0){\line(1,0){3}}
	\put(0,1){\line(1,0){3}}
	\put(3,0){\line(0,1){1}}
	\put(2,1){\line(0,1){1}}
	\put(1,2){\line(1,0){1}}
	\put(0,3){\line(1,0){1}}
	\put(1,2){\line(0,1){1}}
	\put(0,0){\line(0,1){3}}
	\end{picture}%
	\qquad
	\begin{picture}(3,3)
	\put(0,0){\line(1,0){3}}
	\put(0,1){\line(1,0){1}}
	\put(2,1){\line(1,0){1}}
	\put(3,0){\line(0,1){1}}
	\put(1,2){\line(1,0){1}}
	\put(2,1){\line(0,1){1}}
	\put(1,0){\line(0,1){1}}
	\put(0,3){\line(1,0){1}}
	\put(1,2){\line(0,1){1}}
	\put(0,0){\line(0,1){3}}
	\end{picture}%
	\qquad
\end{equation*}
with $\sign(T) \weight(T)=(3,2,1),-(1,4,1),-(3,0,3),(1,0,5)$, 
and $\mathrm{type}(T) = (3,2,1), (4,1,1), (3,3),(5,1)$, respectively.
}
\end{example}

E{\u{g}}ecio{\u{g}}lu and Remmel~\cite{ERemmel:1990} proved that
\begin{equation} \label{equation.Kinv}
	K^{-1}_{\alpha \mu} = \sum_T \sign(T),
\end{equation}
where the sum is over all tabloids $T$ of type $\alpha$ and shape $\mu$.

\section{Cyclic symmetry}
\label{section.cyclic}

To compute the fusion coefficients, it suffices to 
calculate $s_\lambda s_\mu$ in the fusion ring.
Using~\eqref{equation.sh} we obtain
\begin{equation} \label{equation.ss}
	s_\lambda s_\mu = \sum_\alpha K_{\alpha\mu}^{-1} \; h_\alpha s_\lambda \;.
\end{equation}
Note that $h_\alpha$ is multiplicative with
$h_\alpha = h_{\alpha_1} h_{\alpha_2 } \cdots$.
This enables us to compute the product $h_\alpha s_\lambda$ using
the fusion Pieri rule~\cite[Proposition 2.6]{GW:1990}:
for $1\le r\le \ell$ and $\lambda\in \Partition^{\ell,n}$
\begin{equation} \label{equation.fusion pieri}
	h_r s_\lambda = \sum_\nu s_\nu\;,
\end{equation}
where the sum is over all $\nu\in \Partition^{\ell,n}$ such that $\nu/\lambda$ is a horizontal $r$-strip and
$\nu_1-\lambda_n \le \ell$.

\subsection{Cylindric tableaux}
This leads us to the definition of \textit{cylindric tableaux}. 
See also~\cite{GK:1997,Korff:2011,McNamara:2006,Postnikov:2005}. For the precise definition we 
use a notion of shifting (skew) partitions. View a skew 
partition $\nu/\lambda$ as being placed at the
origin so that the bottom leftmost cells of $\nu$ and $\lambda$ 
are placed at $(0,0)$. We then define $\shift_{a,b}(\nu/\lambda)$
to be the skew partition where the bottom leftmost cells of $\nu$ 
and $\lambda$ are placed at position $(a,b)$ in the plane.
We denote the superposition of a skew partition $\nu/\lambda$ 
and its shift by $\nu/\lambda + \shift_{a,b}(\nu/\lambda)$.
We can similarly shift skew tableaux, which are just fillings of skew shapes.

Note that when $\nu, \lambda\in \Partition^{\ell,n}$ such that $\nu/\lambda$ 
is a skew shape, then $\nu/\lambda+\shift_{-\ell,n}(\nu/\lambda)$ 
can be viewed as a skew shape inside the quadrant $x\ge -\ell$ and $y\ge 0$.

\begin{definition}
For two positive integers $\ell$ and $n$ and $\lambda\subseteq\nu\in 
\Partition^{\ell,n}$, a cylindric tableau $t$ of shape $\nu/\lambda$ is 
a column-strict filling of the shape $\nu/\lambda$ such that 
$t$+$\shift_{-\ell,n}(t)$ is still column-strict. 

We denote the set of all cylindric tableaux of shape $\nu/\lambda$ and 
content $\mu$ by $\Tab_{\nu/\lambda,\mu}^{\cyc}$ where,
as usual, the content $\mu$ of a tableau $t$ is the tuple such that $\mu_i$ 
is the number of letters $i$ in $t$.
\end{definition}

\begin{example} \label{example.cyclic}
{\rm
Let $\lambda=(1,1)$, $\mu=(2,2,2)$, $\nu=(4,2,2)$ and $\ell=n=3$. 
Then there are two cylindric tableaux
in $\Tab_{\nu/\lambda,\mu}^{\cyc}$:
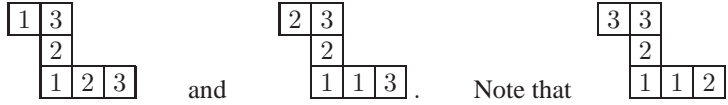
\[
	\tableau[sbY]{1,3|\bl,2|\bl,1,2,3} \qquad \text{and} \qquad
	\tableau[sbY]{2,3|\bl,2|\bl,1,1,3} \; .
	\qquad \text{Note that} \quad
	\tableau[sbY]{3,3|\bl,2|\bl,1,1,2} 
\]
is not cylindric since after shifting by $(-3,3)$ the rightmost 2 would sit above the leftmost 3, which is not column-strict any longer.
}
\end{example}

\subsection{Fusion coefficients}

By iteration of~\eqref{equation.fusion pieri}, we derive that
\begin{equation}
	h_\alpha s_\lambda = \sum_\nu K_{\nu/\lambda, \alpha}^{\cyc} s_\nu \;,
\end{equation}
where $K_{\nu/\lambda,\alpha}^{\cyc} = |\Tab_{\nu/\lambda,\alpha}^\cyc|$ is the cardinality of the set of cylindric tableaux of 
skew shape $\nu/\lambda$ and content $\alpha$. Combining this with~\eqref{equation.ss} we obtain
\begin{equation}
	s_\lambda s_\mu = \sum_{\alpha,\nu} K_{\alpha\mu}^{-1} \; K_{\nu/\lambda, \alpha}^{\cyc} s_\nu\; ,
\end{equation}
which shows that the fusion coefficient is given by the formula
\begin{equation} \label{equation.fusion signed}
	c_{\lambda \mu}^{\nu,\ell,n} =\sum_\alpha K_{\nu/\lambda, \alpha}^{\cyc} \; K_{\alpha\mu}^{-1} \;.
\end{equation}
As we can see from~\eqref{equation.Kinv}, the inverse of the Kostka matrix contains negative signs, so this 
formula is an alternating sum. In the next section we will discuss a sign-reversing involution to cancel
terms in certain cases; the number of fixed points under this involution will precisely amount to the fusion coefficient.

\begin{example} \label{example.fusion}
{\rm 
As in Example~\ref{example.cyclic} consider $\lambda=(1,1)$, $\mu=(2,2,2)$, $\nu=(4,2,2)$ and $\ell=n=3$.
The tabloids of shape $\mu$ are
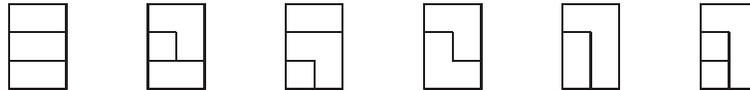
\[
	\unitlength=\cellsize
	\begin{picture}(3,3)
	\put(0,0){\line(1,0){2}}
	\put(0,1){\line(1,0){2}}
	\put(0,2){\line(1,0){2}}
	\put(0,3){\line(1,0){2}}
	\put(0,0){\line(0,1){3}}
	\put(2,0){\line(0,1){3}}
	\end{picture}%
	\qquad
	\begin{picture}(3,3)
	\put(0,0){\line(1,0){2}}
	\put(0,1){\line(1,0){2}}
	\put(0,2){\line(1,0){1}}
	\put(0,3){\line(1,0){2}}
	\put(0,0){\line(0,1){3}}
	\put(2,0){\line(0,1){3}}
	\put(1,1){\line(0,1){1}}
	\end{picture}%
	\qquad
	\begin{picture}(3,3)
	\put(0,0){\line(1,0){2}}
	\put(0,1){\line(1,0){1}}
	\put(0,2){\line(1,0){2}}
	\put(0,3){\line(1,0){2}}
	\put(0,0){\line(0,1){3}}
	\put(1,0){\line(0,1){1}}
	\put(2,0){\line(0,1){3}}
	\end{picture}%
	\qquad
	\begin{picture}(3,3)
	\put(0,0){\line(1,0){2}}
	\put(1,1){\line(1,0){1}}
	\put(0,2){\line(1,0){1}}
	\put(0,3){\line(1,0){2}}
	\put(0,0){\line(0,1){3}}
	\put(1,1){\line(0,1){1}}
	\put(2,0){\line(0,1){3}}
	\end{picture}%
	\qquad
	\begin{picture}(3,3)
	\put(0,0){\line(1,0){2}}
	\put(0,2){\line(1,0){1}}
	\put(0,3){\line(1,0){2}}
	\put(0,0){\line(0,1){3}}
	\put(1,0){\line(0,1){2}}
	\put(2,0){\line(0,1){3}}
	\end{picture}%
	\qquad
	\begin{picture}(3,3)
	\put(0,0){\line(1,0){2}}
	\put(0,1){\line(1,0){1}}
	\put(0,2){\line(1,0){1}}
	\put(0,3){\line(1,0){2}}
	\put(0,0){\line(0,1){3}}
	\put(1,0){\line(0,1){2}}
	\put(2,0){\line(0,1){3}}
	\end{picture}%
\]
with $\sign(T)\weight(T)=(2,2,2),-(2,1,3),-(1,3,2),(0,3,3),-(0,2,4),(1,1,4)$, respectively.
There are no cylindric tableaux of skew shape $\nu/\lambda$ and weight $(0,3,3),(0,2,4)$ or $(1,1,4)$.
The cylindric tableaux of skew shape $\nu/\lambda$ and weights $(2,2,2),(2,1,3)$, and $(1,3,2)$ are
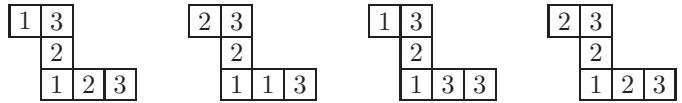
\[
	\tableau[sbY]{1,3|\bl,2|\bl,1,2,3} \qquad
	\tableau[sbY]{2,3|\bl,2|\bl,1,1,3} \qquad
	\tableau[sbY]{1,3|\bl,2|\bl,1,3,3} \qquad
	\tableau[sbY]{2,3|\bl,2|\bl,1,2,3} \; .
\]
Since two of them come with a positive sign and two with a negative sign, Equation~\eqref{equation.fusion signed}
shows that the fusion coefficient $c_{(11),(222)}^{(422),3,3} =0$.
}
\end{example}

Note that if instead of the fusion Pieri rule as 
in~\eqref{equation.fusion pieri} one uses the usual Pieri rule for $h_r s_\lambda$,
one obtains the following expression for the Littlewood--Richardson coefficients
\begin{equation} \label{equation.LR signed}
	c_{\lambda \mu}^{\nu} =\sum_\alpha K_{\nu/\lambda, \alpha} \; K_{\alpha\mu}^{-1} \;,
\end{equation}
where $K_{\nu/\lambda, \alpha}$ is the skew Kostka matrix.

\section{Sign-reversing involution}
\label{section.remmel}

Remmel and Shimozono~\cite{RS:1998} proved the Littlewood--Richardson rule using a sign-reversing
involution. Let us explain their approach first as we will use a modification of it in Section~\ref{subsection.proof}
to prove Theorem~\ref{theorem.fusion}.

\subsection{The Littlewood--Richardson case} \label{subsection.LR}
Note that each tabloid $T$ as defined in Section~\ref{section.tabloids} 
is in one-to-one correspondence with its weight.  Clearly a given $T$ 
yields $\weight(T)$. Conversely, given a weight 
$(\alpha_1,\alpha_2,\ldots)$, start with
the bottommost cell in the first column and draw a ribbon of length 
$\alpha_1$; there is a unique way of doing so.
Then proceed to the second cell in the first column and 
draw a ribbon of length $\alpha_2$ and so on. It is not
hard to see that either there is no way of doing so or there is
a unique way of drawing the ribbons such that
at each step the resulting shape is a partition. 

Under this correspondence between tabloids and weights, \eqref{equation.LR signed} can be rewritten as
(see also~\cite[Eq. (1.14)]{RS:1998})
\begin{equation} \label{equation.c alt}
	c_{\lambda \mu}^\nu = \sum_{(\sigma,t)} \sign(\sigma) \; ,
\end{equation}
where the sum is over all pairs $(\sigma,t)$ with $\sigma \in S_n$ and $t\in \Tab_{\nu/\lambda,\alpha}$ a column-strict
skew tableau of shape $\nu/\lambda$ and weight $\alpha=\sigma(\mu+\rho)-\rho$. Here 
$\rho=(n-1,n-2,\ldots,1,0)$ and $\sigma$ acts on an $n$-tuple by permuting its entries.

The set of column-strict skew tableaux of given shape $\nu/\lambda$ 
over the alphabet $\{1,2,\ldots,n\}$, denoted $\Tab_{\nu/\lambda}$, 
is endowed with \textit{crystal operators} $\et_i,\ft_i,\st_i$ for $1\le i<n$.
For $t\in \Tab_{\nu/\lambda}$, let $\word(t)$ be the column reading 
word of $t$. That is, read the columns of $t$ top to bottom, left 
to right. On a word $w$, the crystal operators $\et_i,\ft_i$, and $\st_i$ 
only act on the letters $i$ and $i+1$. In the subword of $w$ consisting of the 
letters $i$ and $i+1$, successively bracket pairs
$(i+1) \; i$. Then $\ft_i$ makes the rightmost unbracketed $i$ into an 
$i+1$; if no such $i$ exists, $\ft_i$ annihilates the word. Similarly, 
$\et_i$ changes the leftmost unbracketed $i+1$ into an $i$; if no 
such $i+1$ exists, $\et_i$ annihilates the word. Finally, if the subword 
of $w$ of unbracketed letters $i$ and $i+1$ is $i^a (i+1)^b$, then in $\st_i(w)$
this subword is replaced by $i^b (i+1)^a$.

\begin{example}
{\rm
Let $n=4$ and $w=4123322341214223$. The bracketing for $i=2$ yields
\begin{equation*}
\begin{split}
	& \begin{array}{cccccccccccccccc}
	4&1&2      &3&3&2&2&3&4&1&2&1&4& 2       & 2 & 3\\
	  &  &\cdot&( &( & )& )& (&   &  & )&  &  & \cdot & \cdot & \cdot
	\end{array}\\
	\intertext{so that}
	\qquad \qquad \et_2(w) = 
	& \begin{array}{cccccccccccccccc}
	4&1&2      &3&3&2&2&3&4&1&2&1&4& 2       & 2&2
	\end{array}\\
	\ft_2(w) = 
	& \begin{array}{cccccccccccccccc}
	4&1&2      &3&3&2&2&3&4&1&2&1&4& 2       & 3& 3
	\end{array}\\
	\st_2(w) = 
	& \begin{array}{cccccccccccccccc}
	4&1&2      &3&3&2&2&3&4&1&2&1&4& 3       & 3& 3
	\end{array}\; .
\end{split}
\end{equation*}	
}
\end{example}
The action of the crystal operators on column-strict skew tableaux is determined by the action on their column-words.
It is known that after the application of $\et_i,\ft_i$, and $\st_i$, the resulting skew tableau is still column-strict.
A tableau $t\in \Tab_{\nu/\lambda}$ is called \textit{highest weight} if $\et_i(t)=0$ for all $1\le i<n$.
Notice that $t$ is highest weight if and only if $\word(t)$ is lattice.

To prove the Littlewood--Richardson rule, Remmel and Shimozono~\cite{RS:1998} introduced the following sign-reversing
involution $\theta$:
\begin{enumerate}
\item If $t$ is highest weight, then $\theta(\sigma,t)=(\sigma,t)$.
\item Otherwise, let $r+1$ be the rightmost letter in $\word(t)$ that violates the 
lattice condition. Define
$\theta(\sigma,t) = (\sigma_r \sigma, \st_r \et_r(t))$.
\end{enumerate}
Since the highest weight elements $t\in \Tab_{\nu/\lambda}$ yield the fixed-points of $\theta$, their weight must be
a partition. It was shown in~\cite{RS:1998}, that this implies that $\sigma=\id$ so that $t\in \Tab_{\nu/\lambda,\mu}$.
Hence indeed $c_{\lambda \mu}^\nu$ counts the tableaux in $\Tab_{\nu/\lambda,\mu}$ that are lattice (or equivalently
highest weight).

\subsection{Proof of Theorem~\ref{theorem.fusion}} \label{subsection.proof}
Before giving the proof of Theorem~\ref{theorem.fusion}, we need to provide the definition of \textit{cutting points} and
$\tilde{\nu}$ and $\tilde{\lambda}$ in the statement of the theorem.

\begin{definition} \label{definition.cut}
Given two positive integers $\ell$ and $n$, 
consider $\lambda\in\Rec^{\ell,n}$ and $\nu\in\Partition^{\ell,n}$ 
such that $\nu/\lambda$ is a skew shape where
$\nu/\lambda + \shift_{-\ell,n}(\nu/\lambda)$ is not connected.
A cutting point $c$ is the index of the rightmost column 
in $\nu/\lambda + \shift_{-\ell,n}(\nu/\lambda)$ such that the columns with 
$x$-coordinate
$c$ and $c+1$ do not share a common edge.  Note that such a $c$
must exist since the skew shape is not connected.  Let $\tilde{\nu}$ and 
$\tilde{\lambda}$ denote the partitions where
$\tilde{\nu}/\tilde{\lambda}$ is the
skew shape in the window from column $c-\ell$ to column $c$.
\end{definition}

\begin{example}
{\rm
Let $\ell=4$, $n=3$, $\lambda=(3,1)$, and $\nu=(5,5,1)$, so that
\[
	\raisebox{5mm}{$\nu/\lambda =$} \; 
	{\small \tableau[sbY]{|\bl,,,,|\bl,\bl,\bl,,}} 
	\qquad \raisebox{5mm}{and $\qquad \nu/\lambda + \shift_{-4,3}(\nu/\lambda)=$}
	{\small \tableau[sbY]{|\bl,,,,|\bl,\bl,\bl,,|\bl,\bl,\bl,\bl,|\bl,\bl,\bl,\bl,\bl,,,,|\bl,\bl,\bl,\bl,\bl,\bl,\bl,,}} 
	\; \raisebox{5mm}{.}
\]
Then $c=1$ is a cutting point and 
\[
	\raisebox{5mm}{$\tilde{\nu}/\tilde{\lambda}=$} \;
	{\small \tableau[sbY]{,,,|\bl,\bl,,|\bl,\bl,\bl,}}
	\; \raisebox{5mm}{.}
\]
Indeed with $\mu=(4,2,1)$ we have that $c_{(31),(421)}^{(551),4,3}=c_{(32),(421)}^{(444)}=1$ verifying 
Theorem~\ref{theorem.fusion}~\eqref{thm.cutting}.
}
\end{example}

For the proof of Theorem~\ref{theorem.fusion} \eqref{thm.cutting}, one may use the same arguments as in the 
derivation of~\eqref{equation.c alt} for the Littlewood--Richardson coefficients to rewrite~\eqref{equation.fusion signed} 
for the fusion coefficients as
\begin{equation}
	c_{\lambda \mu}^{\nu,\ell,n} = \sum_{(\sigma,t)} \sign(\sigma) \;,
\end{equation}
where now the sum is over all pairs $(\sigma,t)$ with $\sigma\in S_n$ and $t\in \Tab_{\nu/\lambda,\alpha}^\cyc$
a cylindric tableau of shape $\nu/\lambda$ and weight $\alpha=\sigma(\mu+\rho)-\rho$.

Due to the cylindric symmetry of the tableaux in $\Tab_{\nu/\lambda,\alpha}^\cyc$, 
we have $|\Tab_{\nu/\lambda,\alpha}^\cyc| = |\Tab_{\tilde{\nu}/\tilde{\lambda},\alpha}^\cyc|$, 
where recall the definition of $\tilde{\nu}$ and $\tilde{\lambda}$ from Definition~\ref{definition.cut}. 
Since $c$ (as in Definition~\ref{definition.cut})
is a cutting point, that is, the adjacent columns do not share an edge, there 
is no cylindric column-strict condition imposed
on the elements in $\Tab_{\tilde{\nu}/\tilde{\lambda},\alpha}^\cyc$. 
Hence $\Tab_{\tilde{\nu}/\tilde{\lambda},\alpha}^\cyc
=\Tab_{\tilde{\nu}/\tilde{\lambda},\alpha}$, so that the arguments from Section~\ref{subsection.LR} apply.
This proves Theorem~\ref{theorem.fusion}~\eqref{thm.cutting}.

\begin{definition} \label{definition.sl2}
Given two positive integers $\ell$ and $n\ge 2$, 
let $\lambda,\mu\in\Rec^{\ell,n}$ and $\nu\in\Partition^{\ell,n}$ be
such that $\mu$ has at most two parts. Then either $\nu/\lambda$ has a cutting point or it contains
at least one column of height two. Let $c$ be the rightmost such column. Then $\tilde{\nu}$ and $\tilde{\lambda}$
denote the partitions such that $\tilde{\nu}/\tilde{\lambda}$ is the skew shape of $\nu/\lambda + \shift_{-\ell,n}(\nu/\lambda)$ 
in the window from column $c-\ell$ to column $c$.
\end{definition}

By the same arguments as in the proof of part~\eqref{thm.cutting} of Theorem~\ref{theorem.fusion} we have
$|\Tab_{\nu/\lambda,\alpha}^\cyc| = |\Tab_{\tilde{\nu}/\tilde{\lambda},\alpha}^\cyc|$ by cyclic symmetry.
Since $\mu$ has only two parts, the only crystal operators that apply in this case are $\st_1 \et_1$. The column
of height 2 (which after the cyclic shift is the rightmost column) contains the letters $21$ by column-strictness.
The crystal operators cannot change this column due to the crystal bracketing rules. Hence the cylindric column-strict
conditions are always guaranteed which proves Theorem~\ref{theorem.fusion}~\eqref{thm.sl2}.
Note that this case is related to the cyclic $\mathfrak{sl}_2$ crystals introduced in~\cite{HS:2012}.

\section{Beyond the cutting point}
\label{section.beyond}

The cylindric tableaux of Example~\ref{example.fusion} do not have a cutting point and hence Theorem~\ref{theorem.fusion}
does not apply. Under the Remmel--Shimozono involution we have
\[
	\raisebox{3mm}{$\st_2\et_2$} \quad  \tableau[sbY]{1,3|\bl,2|\bl,1,2,3} 
	\quad \raisebox{3mm}{$=$} \quad \tableau[sbY]{1,3|\bl,2|\bl,1,3,3} \; .
\]
However, the action of the Remmel--Shimozono involution on the other two cylindric tableaux yields non-cylindric tableaux:
\[
	\raisebox{3mm}{$\st_2 \et_2$} \quad \tableau[sbY]{2,3|\bl,2|\bl,1,1,3}
	\quad \raisebox{3mm}{$=$} \quad \tableau[sbY]{3,3|\bl,2|\bl,1,1,3} 
	\qquad \raisebox{3mm}{and} \qquad
	\raisebox{3mm}{$\st_2 \et_2$} \quad \tableau[sbY]{2,3|\bl,2|\bl,1,2,3}
	\quad \raisebox{3mm}{$=$} \quad \tableau[sbY]{3,3|\bl,2|\bl,1,3,3} 
\]
and hence does not yield a cancelation within the set $\Tab_{\nu/\lambda,\cdot}^\cyc$.

Note, however, that it is possible to amend the operators used by Remmel and Shimozono by conjugating
the action of $\st_i \et_i$ by a cyclic shift. In the above example, moving the 2 in the leftmost column down we 
obtain
\[
	\tableau[sbY]{3|2|1,1,3|\bl,\bl,2} 
	\qquad \raisebox{3mm}{and} \qquad
	\tableau[sbY]{3|2|1,2,3|\bl,\bl,2}
\]
which cancel under the action of $\st_1 \et_1$.

We conjecture that such cyclic cancelations are always possible. In fact, 
computer experiments using {\tt Sage}~\cite{sage,sagecombinat} suggest 
that the resulting fusion lattice tableaux (that is, the skew 
tableaux that are fixed points under the involution) correspond 
to the 2d puzzles conjectured by 
Knutson~\cite{BKT:2003} to yield the quantum Littlewood--Richardson 
coefficients or equivalently fusion 
coefficients by a bijection between puzzles and tableaux 
similar to~\cite[Figure 11]{Vakil:2006}.	


\end{document}